\documentstyle[11pt,amstex,amssymb]{amsart}\textheight=19.8truecm\textwidth=13.7truecm
\begin{document} 

\vskip 20pt
\centerline{\Large\bf  Convexity or concavity inequalities for}
\vskip 20pt
\centerline{\Large\bf  Hermitian operators }
\vskip 40pt
\centerline{Jean-Christophe Bourin}
\vskip 15pt
\centerline{Les Coteaux, 8 rue Henri Durel 78510 Triel, France}
\vskip 10pt
 \centerline{E-mail: bourinjc@@club-internet.fr}
\vskip 25pt
\noindent
{\small {\bf Abstract.} Given a  Hermitian operator, a monotone convex function $f$ and a subspace ${\cal E}$, $\dim{\cal E}<\infty$, there exists a unitary operator $U$ on ${\cal E}$ such that $f(A_{\cal E}) \le Uf(A)_{\cal E}U^*$. (Here $X_{\cal E}$ denotes the compression of $X$ onto ${\cal E}$). A related result is: For a monotone convex function $f$, $0<\alpha,\beta<1$, $\alpha+\beta=1$, and Hermitian operators $A$, $B$ on a finite dimensional space, there exists a unitary $U$ such that $f(\alpha A+\beta B)\le U\{\alpha f(A)+ \beta f(B)\}U^*$. More general convexity results are established. Also, several old and new trace inequalities of Brown-Kosaki and Hansen-Pedersen type are derived. We study the behaviour of the map $p\longrightarrow \{(A^p)_{\cal E}\}^{1/p}$, $A\ge0$, $0<p<\infty$. 


Mathematical subjects classification: 47A20, 47A30, 47A63 }

\vskip 25pt
{\large\bf \ Introduction}
\vskip 10pt
Given an operator $A$ on a separable Hilbert space ${\cal H}$ and a subspace ${\cal E}\subset{\cal H}$, we denote by $A_{\cal E}$ the compression of $A$ onto  ${\cal E}$, i.e.\  the restriction of $EAE$
to ${\cal E}$, $E$ being the projection onto ${\cal E}$. If ${\cal E}$ is a finite dimensional  subspace, we  show  that, for any  Hermitian operator $A$ and any monotone convex function $f$ defined on the spectrum of $A$, there exits a unitary operator $U$ on ${\cal E}$ such that the operator inequality
\vskip 5pt
\begin{equation*}
f(A_{\cal E}) \le Uf(A)_{\cal E}U^*. \tag{$*$}
\end{equation*}
\vskip 5pt\noindent
holds. Here, $f(A)_{\cal E}$ must be read as $(f(A))_{\cal E}$. This result together with the elementary method of its proof  motivate the  whole paper. 
In Section 1 we prove the above inequality and give a version. We also study the map $p\longrightarrow \{(A^p)_{\cal E}\}^{1/p}$, $0<p<\infty$ for a positive operator $A$ on a finite dimensional space and $d$-dimensional subspace ${\cal E}$. In general, this map converges to an operator $B$ on ${\cal E}$ whose eigenvalues are the $d$ largest eigenvalues of $A$.

\newpage
Section 2 is concerned with eigenvalues inequalities (equivalently operator inequalities) which improve some trace inequalities of Brown-Kosaki and Hansen-Pedersen: Given a monotone convex function $f$   defined on the real line with $f(0)\le 0$, a Hermitian operator $A$ and a contractive operator $Z$ acting on a finite dimensional space, there exists a unitary operator $U$ such that
$$
f(Z^*AZ)\le UZ^*f(A)ZU^*.
$$

 In Section 3, we prove that
$$
{\rm Tr\,}f(Z^*AZ)\le {\rm Tr\,}Z^*f(A)Z
$$
for every positive operator $A$ and expansive operator $Z$ on a finite dimensional space, and every concave
function $f$ defined on an interval $[0,b]$, $b\ge\Vert Z^*AZ\Vert_{\infty}$, with $f(0)\ge0$ ($\Vert\cdot\Vert_{\infty}$ denotes the usual operator norm). Under the additional assumption that
$f$ is nondecreasing and nonnegative, this trace inequality entails
$$
\Vert f(Z^*AZ)\Vert_{\infty}\le \Vert Z^*f(A)Z\Vert_{\infty}
$$
In a forthcomming project we will extend this result to the infinite dimensional setting. We will also give a version of $(*)$ for infine dimensional subspace ${\cal E}$, by adding  a $rI$ term in the right hand side, where $I$ stands for the identity and $r>0$ is arbitrarily small.

\vskip 20pt \noindent
{\large\bf 1. Compressions and convex functions}
\vskip 10pt

By a classical result of C.\ Davis [4] (see also [1, p.\ 117-9]), a function $f$ on $(a,b)$ is operator convex if and only if for every subspace ${\cal E}$ and every Hermitian operator $A$ whose spectrum lies in  $(a,b)$  one has
\begin{equation}
f(A_{\cal E}) \le f(A)_{\cal E}. \tag 1
\end{equation} 
What can be said about convex, not operator convex functions ? Let $g$ be  operator convex on $(a,b)$ and let $\phi$ be a nondecreasing, convex function on $g((a,b))$. Then, $f=\phi\circ g$ is convex and we say that $f$ is {\it unitary} convex on $(a,b)$. Since $t\longrightarrow -t$ is trivially operator convex, we note that {\it the class of unitary convex functions contains the class of monotone convex functions}. The following result holds:

\vskip 10pt \noindent
{\bf Theorem 1.1.} {\it Let $f$ be a monotone convex, or more generally unitary convex, function on $(a,b)$ and let $A$ be a  Hermitian operator whose spectrum  lies in $(a,b)$. 
 If  ${\cal E}$ is a  finite dimensional subspace, then there exists a unitary operator $U$ on ${\cal E}$ such that
\begin{equation*}
f(A_{\cal E}) \le Uf(A)_{\cal E}U^*.
\end{equation*}
}

\vskip 5pt\noindent
{\bf Proof.} We begin by assuming that $f$ is monotone. Let $d=\dim{\cal E}$ and let $\{\lambda_k(X)\}_{k=1}^d$ denote the eigenvalues of the Hermitian operator $X$ on ${\cal E}$, arranged in decreasing order and counted with their multiplicities. Let $k$ be an integer, $1\le k\le d$.  There exists a spectral subspace ${\cal F}\subset{\cal E}$ for $A_{\cal E}$ (hence for $f(A_{\cal E})$) , $\dim {\cal F}=k$, such that
\begin{align*} \lambda_k[f(A_{\cal E})] &=\min_{h\in{\cal F};\ \Vert h\Vert=1} \langle h,f(A_{\cal F})h \rangle  \\
&= \min\{f(\lambda_1(A_{\cal F}))\,;\,f(\lambda_k(A_{\cal F}))\} \\
&= \min_{h\in{\cal F};\ \Vert h\Vert=1} f(\langle h,A_{\cal F}h \rangle)  \\
&= \min_{h\in{\cal F};\ \Vert h\Vert=1} f(\langle h,Ah \rangle) 
\end{align*}
where at the second and third steps we use the monotony of $f$. The convexity of $f$ implies 
$$
f(\langle h,Ah \rangle) \le \langle h,f(A)h \rangle
$$
for all normalized vectors $h$. Therefore, by the minmax principle,
\begin{align*}
\lambda_k[f(A_{\cal E})] &\le \min_{h\in{\cal F};\ \Vert h\Vert=1} \langle h,f(A)h \rangle \\
&\le \lambda_k[f(A)_{\cal E}].
\end{align*} 
This  statement is equivalent to the existence of a unitary operator $U$ on ${\cal E}$ satisfying  the conclusion of the theorem. 

If $f$ is unitary convex, $f=\phi\circ g$ with $g$ operator convex and $\phi$ nondecreasing convex; inequality (1) applied to $g$ combined with the fact that $\phi$ is nondecreasing yield a unitary operator $V$ on
${\cal E}$ for which
\begin{equation*}
\phi\circ g(A_{\cal E}) \le V\phi[g(A)_{\cal E}]V^*.
\end{equation*} 
Applying  the first part of the proof to $\phi$ gives a unitary operator $W$ on ${\cal E}$ such that
\begin{equation*}
\phi[g(A)_{\cal E}] \le W[\phi\circ g(A)]_{\cal E}W^*.
\end{equation*} 
We then get the result by letting $U=VW$. \qquad $\Box$ 

\vskip 10pt
Later, we will see that Theorem 1.1 can {\it not} be extended to all convex functions $f$ (Example 2.4).  Of course Theorem 1.1 holds with a reverse inequality for monotone concave functions $f$ (or $f=\phi\circ g$, $g$ operator convex and $\phi$ decreasing concave).

\vskip 10pt 
Given a  positive operator $A$ on a finite dimensional space and a subspace ${\cal E}$, it is natural to study the behaviour of the map
$$
p\longrightarrow \{(A^p)_{\cal E}\}^{1/p}
$$
on $(0,\infty)$. 
The notation $A=\sum_k\lambda_k(A)\,f_k \otimes f_k$ means that $f_k$ is a norm one eigenvector associated to $\lambda_k(A)$ and $f_k \otimes f_k$ is the corresponding norm one projection

\newpage

\vskip 10pt\noindent
{\bf Theorem 1.2.} {\it Let $A=\sum_k\lambda_k(A)\,f_k \otimes f_k$ be a positive operator on a finite dimensional space and let ${\cal E}$ be a  subspace.
Assume  ${\cal E}\cap{\rm span}\{f_j:\,j>d\}=0$. Then, for every integer $k\le \dim{\cal E}$,
 the map $p\longrightarrow \lambda_k(\{(A^p)_{\cal E}\}^{1/p})$ increases on $(0,\infty)$ and
\begin{equation*}
\lim_{p\rightarrow \infty} \lambda_k( \{(A^p)_{\cal E}\}^{1/p}) = \lambda_k(A).
\end{equation*} 
Moreover the family $\{(A^p)_{\cal E}\}^{1/p}$ converges in norm when $p\rightarrow\infty$ and the map $p\longrightarrow \{(A^p)_{\cal E}\}^{1/p}$ is increasing for the Loewner order on $[1,\infty)$. }

\vskip 10pt\noindent
{\bf Proof.} Let $p>0$ and $r>1$.  By Theorem 1.1, there exists a unitary $U:{\cal E}\longrightarrow{\cal E}$ such that
$$
\{(A^p)_{\cal E}\}^r \le U(A^{pr})_{\cal E}U^*,
$$
hence, for all $k$,
$$
\lambda_k^r((A^p)_{\cal E}) \le \lambda_k((A^{pr})_{\cal E}),
$$
so,
$$
\lambda_k(\{(A^p)_{\cal E}\}^{1/p}) \le \lambda_k(\{(A^{pr})_{\cal E}\}^{1/pr}),
$$
that is, the map $p\longrightarrow \lambda_k(\{(A^p)_{\cal E}\}^{1/p})$ increases on $(0,\infty)$. In order to study its convergence when $p\to\infty$, we first show that
\begin{equation*}
\lim_{p\rightarrow \infty} \lambda_1((EA^pE)^{1/p}) = \lambda_1(A) \tag{2}
\end{equation*} 
where $E$ denotes the projection onto ${\cal E}$. 
We note that
\begin{equation*}
\lim_{p\rightarrow \infty} \lambda_1((EA^pE)^{1/p}) \le \lambda_1(A). \tag{3}
\end{equation*}
Recall that $A=\sum_k \lambda_k(A)\,f_k\otimes f_k$.  Since by assumption $f_1\not\in{\cal E}^{\perp}$, there exists a normalized vector $g$  in ${\cal E}$ such that $\langle g, f_1\rangle\neq0$.  Setting $G= g\otimes g$, we have
\begin{equation*}
\lambda_1((GA^pG)^{1/p}) = \langle g, A^pg\rangle^{1/p}= \bigl(\,
\sum_k \lambda_k^p(A)|\langle g, f_k\rangle|^2\,\bigr)^{1/p}.
\end{equation*} 
The above expression is a weighted $l^p$-norm of the sequence $\{ \lambda_k(A)\}$. When $p\rightarrow\infty$, this tends towards the  $l^{\infty}$-norm which is $\lambda_1(A)$. Since
$$
\lambda_1((GA^pG)^{1/p}) \le \lambda_1((EA^pE)^{1/p})
$$
we then deduce with (3) that (2) holds.

 In order to prove the general limit assertion, we consider antisymmetric tensor products. Let $F$ be the projection onto ${\cal F}={\rm span}\{f_j:\, j\le \dim{\cal E}\}$. By assumption $F$ maps ${\cal E}$ onto ${\cal F}$. Therefore $\wedge^k(F)$ maps $\wedge^k({\cal E})$ onto $\wedge^k({\cal F})$ and 
 we may find a norm one tensor $\gamma\in\wedge^k({\cal E})$ such that $\langle\gamma,f_1\wedge\dots\wedge f_k\rangle\neq0$. Hence, with
$\wedge^kE$ and $\wedge^kA$ in place of $E$ and $A$, $1\le k\le \dim{\cal E}$, we may apply (2) to obtain
\begin{equation*}
\lim_{p\rightarrow \infty} \lambda_1(\wedge^k(EA^pE)^{1/p}) = \lambda_1(\wedge^kA)
\end{equation*}  
meaning that
\begin{equation*}
\lim_{p\rightarrow \infty} \prod_{1\le j\le k} \lambda_j((EA^pE)^{1/p}) = \prod_{1\le j\le k}\lambda_j(A). 
\end{equation*}  
From these relations we infer that, for every $k\le \dim{\cal E}$, we have
\begin{equation*}
\lim_{p\rightarrow \infty} \lambda_k((EA^pE)^{1/p}) = \lambda_k(A) 
\end{equation*} 
proving the main assertion of the theorem.

For $p,r\ge 1$ we have
$$
(EA^{pr}E)^{1/r}\ge EA^pE.
$$
by Hansen's inequality [6]. Since $t\longrightarrow t^{1/p}$ is operator monotone by the
Loewner theorem [9, p.\ 2], we have
$$
(EA^{pr}E)^{1/pr}\ge (EA^pE)^{1/p}.
$$
Thus $p\longrightarrow (EA^pE)^{1/p}$ increases on $[1,\infty)$. Since this map is  bounded, it
 converges in norm.  \qquad $\Box$
 
 \vskip 15pt
 The author is indebted to a referee for having pointed out a misconception in the initial statement and proof of Theorem 1.2.

\vskip 20pt \noindent
{\large\bf 2. Contractions and convex functions}

 \vskip 10pt
In [6] and [7], the authors show that inequality (1) is equivalent to the following statement.

\vskip 10pt \noindent
{\bf Theorem 2.1.} (Hansen-Pedersen) {\it Let $A$ 	and  $\{A_i\}_{i=1}^m$ be  Hermitian operators  and  let $f$ be an operator convex function
 defined on an interval $[a,b]$ containing the spectra of $A$ and $A_i$, $i=1,\dots m$.
 \vskip 5pt
 $(1)$ If $Z$ is a contraction, $0\in[a,b]$  and $f(0)\le 0$, 
\begin{equation*}
 f(Z^*AZ) \le Z^*f(A)Z.
\end{equation*}
\vskip 5pt
 $(2)$ If $\{Z_i\}_{i=1}^m$ is an isometric column, 
\begin{equation*}
 f(\sum_iZ_i^*A_iZ_i) \le \sum_iZ_i^*f(A_i)Z_i.
\end{equation*}
}
\vskip 5pt \noindent
Here, an isometric column $\{Z_i\}_{i=1}^m$ means that $\sum_{i=1}^mZ_i^*Z_i=I$. 

\vskip 10pt
In a similar way, Theorem 1.1 is equivalent to the next one.
We state it in the finite dimensional setting, but an analogous version exists in the  infinite dimensional  setting by adding a $rI$ term in the right hand side of the inequalities.

\vskip 10pt \noindent
{\bf Theorem 2.2.} {\it Let $A$ 	and  $\{A_i\}_{i=1}^m$ be  Hermitian operators on a finite dimensional space and  let $f$ be a monotone, or more generally unitary, convex function
 defined on an interval $[a,b]$ containing the spectra of $A$ and $A_i$, $i=1,\dots m$.
 \vskip 5pt
 $(1)$ If $Z$ is a contraction, $0\in[a,b]$  and $f(0)\le 0$, then there exists a unitary operator $U$ such that
\begin{equation*}
 f(Z^*AZ) \le UZ^*f(A)ZU^*.
\end{equation*}
\vskip 5pt
 $(2)$ If $\{Z_i\}_{i=1}^m$ is an isometric column, then there exists a unitary operator $U$ such that
\begin{equation*}
 f(\sum_iZ_i^*A_iZ_i) \le U\{\sum_iZ_i^*f(A_i)Z_i\}U^*.
\end{equation*}
}

Here, we give a first proof based on Theorem 1.1. A  more direct proof is given at the end of the section.

\vskip 10pt \noindent
{\bf Proof.} Theorem 2.2 and  Theorem 1.1 are equivalent. Indeed, to prove Theorem 1.1 we may  first assume, by a limit argument, that $f$ is defined on the whole real line. Then, we may assume that $f(0)=0$ so that Theorem 1.1 follows from Theorem 2.2 by taking $Z$ as the projection onto ${\cal E}$.

 Theorem 1.1 entails Theorem 2.2(1): to see that, we introduce the partial isometry $V$ and the operator $\tilde{A}$ on ${\cal H}\oplus{\cal H}$ defined by
$$ V=
\begin{pmatrix}
Z&0 \\ (I-|Z|^2)^{1/2}&0
\end{pmatrix},\quad
\tilde{A}=
\begin{pmatrix}
A&0 \\ 0&0
\end{pmatrix}.
$$
Denoting by ${\cal H}$ the first summand of the direct sum ${\cal H}\oplus{\cal H}$, we observe that
$$
f(Z^*AZ)= f(V^*\tilde{A}V)\!:\!{\cal H} = V^*f(\tilde{A}_{V({\cal H})})V\!:\!{\cal H}.
$$
Applying  Theorem 1.1 with ${\cal E}=V({\cal H})$, we get a unitary operator $W$ on  $V({\cal H})$ such that
$$
f(Z^*AZ) \le V^*Wf(\tilde{A})_{V({\cal H})}W^*V\!:\!{\cal H}. 
$$
Equivalently, there exists a unitary operator $U$ on ${\cal H}$ such that
\begin{align*}
f(Z^*AZ) \le& UV^*f(\tilde{A})_{V({\cal H})}(V\!:\!{\cal H})U^* \\
=&UV^* \begin{pmatrix} f(A)&0 \\0&f(0) \end{pmatrix}(V\!:\!{\cal H})U^* \\
=& U\{Z^*f(A)Z + (I-|Z|^2)^{1/2}f(0)(I-|Z|^2)^{1/2}\}U^*.
\end{align*}
Using $f(0)\le 0$ we obtain the first claim of Theorem 2.2.

 Similarly, Theorem 1.1 implies Theorem 2.2(2) (we may assume $f(0)=0$) by considering the partial isometry and the operator on $\oplus^m{\cal H}$,
$$
\begin{pmatrix}
Z_1&0&\cdots&0 \\ \vdots &\vdots &\ &\vdots \\
Z_m &0 &\cdots &0
\end{pmatrix},\qquad
\begin{pmatrix}
A_1&\ &\  \\ 
\ &\ddots &\  \\
\ &\  &A_m
\end{pmatrix} .
$$
\qquad $\Box$

\vskip 10pt 
 We note that Theorem 2.2 strengthens some well-known trace inequalities:

\vskip 10pt \noindent
{\bf Corollary  2.3.} {\it Let $A$ 	and  $\{A_i\}_{i=1}^m$ be  Hermitian operators on a finite dimensional space and  let $f$ be a  convex  function
 defined on an interval $[a,b]$ containing the spectra of $A$ and $A_i$, $i=1,\dots m$.
 \vskip 5pt
 $(1)$ {\rm (Brown-Kosaki [2])} If $Z$ is a contraction, $0\in[a,b]$  and $f(0)\le 0$, then 
\begin{equation*}
 {\rm Tr\,}f(Z^*AZ) \le {\rm Tr\,}Z^*f(A)Z.
\end{equation*}
\vskip 5pt
 $(2)$ {\rm (Hansen-Pedersen [7])} If $\{Z_i\}_{i=1}^m$ is an isometric column, then 
\begin{equation*}
 {\rm Tr\,}f(\sum_iZ_i^*A_iZ_i) \le {\rm Tr\,}\{\sum_iZ_i^*f(A_i)Z_i\}.
\end{equation*}
} 

\vskip 5pt \noindent
{\bf Proof.} By a limit argument, we may assume that $f$ is defined on the whole real line and can be written as $f(x)=g(x)-\lambda x$ for some convex monotone function $g$ and some scalar $\lambda$. We then apply Theorem 2.2 to $g$. \qquad $\Box$

\vskip 10pt
A very special case of Theorem 2.2(2) is: {\it Given two Hermitian operators $A$, $B$   and a monotone convex or  unitary convex function $f$ on a suitable interval, there exists a unitary operator $U$ such that
\begin{equation*}
 f(\frac{A+B}{2}) \le U\frac{f(A)+f(B)}{2}U^*.
\end{equation*} 
}
This shows that Theorem 2.2, and consequently Theorem 1.1, can not be valid for all convex functions:

\vskip 10pt\noindent
{\bf Example 2.4.} Theorems 1.1 and 2.2 are not valid for a simple convex function such as $t\longrightarrow |t|$. Indeed, it is well-known that the inequality
\begin{equation*}
|A+B|\le U(|A|+|B|)U^* \tag 5
\end{equation*}
is not always true, even for Hermitians $A$, $B$. We reproduce the counterexample [8, p.\ 1]: Take
$$
A=\begin{pmatrix} 1&1 \\ 1&1\end{pmatrix}, \quad B=\begin{pmatrix} 0&0 \\ 0&-2\end{pmatrix}.
$$
Then, as the two eigenvalues of $|A+B|$ equal to $\sqrt{2}$ while $|A|+|B|$ has an eigenvalue equal to $2-\sqrt{2}$, inequality (5) can not hold.

\vskip 10pt
  In connection with Example 2.4, a famous result (e.g., [1, p.\ 74]) states the existence, for any operators $A$, $B$ on a finite dimensional space, of unitary operators $U$, $V$ such that
\begin{equation*}
|A+B|\le U|A|U^*+V|B|V^*. \tag 6
\end{equation*}
\vskip 10pt
In the case of Hermitians $A$, $B$, the above inequality has the following generalization:

\vskip 10pt \noindent
{\bf Proposition 2.5.} {\it Let $A$, $B$ be hermitian operators on a finite dimensional space and let $f$ be an even convex function on the real line. Then, there exist unitary operators $U$, $V$ such that
$$
f\left(\frac{A+B}{2}\right)\le \frac{Uf(A)U^*+Vf(B)V^*}{2}. 
$$
}
\vskip 5pt \noindent
{\bf Proof.} Since $f(X)=f(|X|)$, inequality (6) and the fact that $f$ is increasing on $[0,\infty)$
give unitary operators $U_0$, $V_0$ such that
$$
f\left(\frac{A+B}{2}\right)\le f\left(\frac{U_0|A|U_0^*+V_0|B|V_0^*}{2}\right).
$$
Since $f$ is monotone convex on $[0,\infty)$, Theorem 2.2 completes the proof. \qquad$\Box$

\vskip 10pt \noindent
{\bf Question 2.6.} Does Proposition 2.5 hold for all convex functions defined on the  whole real line ?

\vskip 10pt
 We close this section by giving a direct and  proof of Theorem 2.2, which is a simple adaptation of the proof of Theorem 1.1.

\vskip 10pt \noindent
{\bf Proof.} We restrict ourselves to the case when $f$ is monotone. We will use the following observation  which follows from the standard Jensen's inequality: for any vector $u$ of norm less than or equal to one, since $f$ is convex and $f(0)\le0$,
$$
f(\langle u, Au \rangle) \le \langle u,f(A)u \rangle.
$$
We begin by proving assertion (1).  We have, for each integer $k$ less than or equal to the dimension of the space, a subspace ${\cal F}$ of dimension $k$ such that
\begin{align*}
\lambda_k[f(Z^*AZ)] &= \min_{h\in{\cal F};\ \Vert h\Vert=1} \langle h,f(Z^*AZ)h \rangle \\
&= \min_{h\in{\cal F};\ \Vert h\Vert=1} f(\langle h,Z^*AZh \rangle)  \\
&= \min_{h\in{\cal F};\ \Vert h\Vert=1} f(\langle Zh,AZh \rangle).  
\end{align*}
where we have used the monotony of $f$. Then, using the above observation  and the minmax principle,
\begin{align*}
\lambda_k[f(Z^*AZ)] &\le \min_{h\in{\cal F};\ \Vert h\Vert=1} \langle Zh,f(A)Zh \rangle \\
&\le \lambda_k[Z^*f(A)Z].
\end{align*} 

We turn to assertion (2). For any integer $k$ less than or equal to the dimension of the space, we have a subspace ${\cal F}$ of dimension $k$ such that
\begin{align*}
\lambda_k[f(\sum Z_i^*A_iZ_i)] &= \min_{h\in{\cal F};\ \Vert h\Vert=1} \langle h,f(\sum Z_i^*A_iZ_i)h \rangle \\
&= \min_{h\in{\cal F};\ \Vert h\Vert=1} f(\langle h,\sum Z_i^*A_iZ_ih \rangle)  \\
&= \min_{h\in{\cal F};\ \Vert h\Vert=1} f(\sum \Vert Z_ih\Vert^2 (\langle Z_ih,A_iZ_ih \rangle/ \Vert Z_ih\Vert^2)) \\
&\le \min_{h\in{\cal F};\ \Vert h\Vert=1} \sum \Vert Z_ih\Vert^2 f(\langle Z_ih,A_iZ_ih \rangle/ \Vert Z_ih\Vert^2) \tag 7\\
&\le \min_{h\in{\cal F};\ \Vert h\Vert=1} \sum  \langle Z_ih,f(A_i)Z_ih \rangle) \tag 8\\
&\le \min_{h\in{\cal F};\ \Vert h\Vert=1} \langle h,\sum Z^*_if(A_i)Z_ih \rangle) \\
&\le \lambda_k[\sum Z_i^*f(A_i)Z_i]
\end{align*}
where we have used in (7) and (8) the convexity of $f$. \qquad$\Box$

\vskip 20pt\noindent
{\large\bf 3. Inequalities involving expansive operators} 
\vskip 10pt
In this section we are in the finite dimensional setting.

\vskip 5pt
For two reals $a$, $z$, with $z>1$, we have $f(za)\ge zf(a)$ for every convex function $f$ with $f(0)\le 0$. In view of Theorem 2.2, one might expect the following result: If $Z$ is an expansive  operator  (i.e.\ $Z^*Z\ge I$), $A$ is a Hermitian operator  and $f$ is a convex function with $f(0)\le0$, then there exists a unitary operator $U$ such that
\begin{equation*}
 f(Z^*AZ) \ge UZ^*f(A)ZU^*. \tag{*}
\end{equation*}
But, as we shall see, this is not always true, even for $A\ge0$ and $f$ nonnegative with $f(0)=0$. Let us first note the following remark: 

\vskip 10pt\noindent
{\bf Remark 3.1.} Let $f:[0,\infty)\longrightarrow[0,\infty)$ be a continuous function with $f(0)=0$. If
$$
{\rm Tr\,}f(Z^*AZ) \le {\rm Tr\,}Z^*f(A)Z
$$
for every positive operator $A$ and every contraction $Z$, then $f$ is convex.

\vskip 10pt 
To check this, it suffices to consider: 
\begin{equation*}
A=\begin{pmatrix} x&0 \\0&y \end{pmatrix} \quad{\rm and}\quad Z=\begin{pmatrix} 1/\sqrt2&0 \\1/\sqrt2&0 \end{pmatrix}
\end{equation*}
where $x$, $y$ are arbitrary nonnegative scalars. Indeed, ${\rm Tr\,}f(Z^*AZ)=f((x+y)/2)$ and ${\rm Tr\,}Z^*f(A)Z=(f(x)+f(y))/2$.

\vskip 10pt We may now state

\vskip 10pt \noindent
{\bf Proposition 3.2.} {\it Let $f:[0,\infty)\longrightarrow[0,\infty)$ be a continuous one to one function with $f(0)=0$ and $f(\infty)=\infty$. Then, the following conditions are equivalent:
\vskip 5pt
 $(1)$ The function $g(t)=1/f(1/t)$ is convex on $[0,\infty)$.
\vskip 5pt
 $(2)$ For every positive operator $A$ and every expansive operator $Z$, there exists a unitary operator $U$ such that
\begin{equation*}
 Z^*f(A)Z \le Uf(Z^*AZ)U^*.
\end{equation*}
}

\vskip 5pt \noindent
{\bf Proof.} We may assume that $A$ is invertible. If $g$ is convex, (note that $g$ is also nondecreasing) then Theorem 2.2 entails that
$$
g(Z^{-1}A^{-1}Z^{-1*}) \le U^*Z^{-1}g(A^{-1})Z^{-1*}U
$$
for some unitary operator $U$. Taking the inverses, since $t\longrightarrow t^{-1}$ is operator decreasing on $(0,\infty)$, this is the same as saying
\begin{equation*}
 Z^*f(A)Z \le Uf(Z^*AZ)U^*.
\end{equation*}
The converse direction follows, again by taking the inverses, from the above remark. \qquad $\Box$

\vskip 10pt
It is not difficult to find convex functions $f:[0,\infty)\longrightarrow[0,\infty)$, with $f(0)=0$ which
do not satisfy to the conditions of Proposition 3.2. So, in general, (*) can not hold. Let us give an explicit simple example. 

\vskip 10pt \noindent
{\bf Example 3.3.} Let $f(t)=t+(t-1)_+$ and
\begin{equation*}
A=\begin{pmatrix} 3/2& 0 \\0& 1/2 \end{pmatrix}, \qquad Z=\begin{pmatrix} 2& 1 \\1& 2 \end{pmatrix}.
\end{equation*}
Then $\lambda_2(f(ZAZ))=0.728..<0.767..=\lambda_2(Zf(A)Z)$. So, (*) does not hold.

\vskip 10pt
In spite of the previous example, we have the following positive result:

\vskip 10pt \noindent
{\bf Lemma 3.4.} {\it Let $A$ be a positive operator, let $Z$ be an expansive operator and $\beta$ be a nonnegative scalar. Then, there exists a unitary operator $U$ such that
\begin{equation*}
 Z^*(A-\beta I)_+Z \le U(Z^*AZ-\beta I)_+U^*.
\end{equation*}
}

\vskip 5pt \noindent
{\bf Proof.} We will use the following simple fact: If $B$ is a positive operator with ${\rm Sp}B\subset\{0\}\cup(x,\infty)$, then we also have ${\rm Sp}Z^*BZ\subset\{0\}\cup(x,\infty)$. Indeed $Z^*BZ$ and $B^{1/2}ZZ^*B^{1/2}$ (which is greater than $B$) have the same spectrum.

Let $P$ be the spectral projection of $A$ corresponding to the eigenvalues strictly greater than $\beta$
and let $A_{\beta}=AP$. Since $t\longrightarrow t_+$ is nondecreasing, there exists a unitary operator $V$ such that 
\begin{equation*}
 (Z^*AZ-\beta I)_+\ge V(Z^*A_{\beta}Z-\beta I)_+V^*
\end{equation*}
Since $Z^*(A-\beta I)_+Z=Z^*(A_{\beta}-\beta I)_+Z$ we may then assume that $A=A_{\beta}$. Now, the above simple fact implies
\begin{equation*}
 (Z^*A_{\beta}Z-\beta I)_+=Z^*A_{\beta}Z-\beta Q
\end{equation*}
where $Q={\rm supp}Z^*A_{\beta}Z$ is the support projection of $Z^*A_{\beta}Z$. Hence, it suffices to show the existence of a unitary operator $W$ such that
\begin{equation*}
 Z^*A_{\beta}Z-\beta Q  \ge WZ^*(A_{\beta}-\beta P)ZW^*=WZ^*A_{\beta}ZW^*-\beta WZ^*PZW^*.
\end{equation*}
But, here we can take $W=I$. Indeed, we have
\begin{equation*}
 {\rm supp}Z^*PZ=Q \ (*) \quad{\rm and}\quad {\rm Sp}Z^*PZ\subset\{0\}\cup[1,\infty) \ (**)
\end{equation*}
where ($\ast\ast$) follows from the above simple fact and the identity ($\ast$) from the  observation below with $X=P$ and $Y=A_{\beta}$.
\vskip 5pt\noindent
{\it Observation.}\, If $X$, $Y$ are two positive operators with ${\rm supp}X={\rm supp}Y$, then for every  operator $Z$ we also have ${\rm supp}Z^*XZ={\rm supp}Z^*YZ$. 
\vskip 5pt\noindent
To check this, we establish the corresponding equality for the kernels,
\begin{equation*}
\ker Z^*XZ=\{h\ :\ Zh\in \ker X^{1/2}\} = \{h\ :\ Zh\in \ker Y^{1/2}\}= \ker Z^*YZ.
\end{equation*}
\qquad $\Box$

\vskip 10pt \noindent
{\bf Theorem 3.5.} {\it Let $A$ be a positive operator and $Z$ be an expansive operator. Assume that
$f$ is a continuous function defined on $[0,b]$, $b\ge \Vert Z^*AZ\Vert_{\infty}$. Then,
\vskip 5pt {\rm (1)} If $f$ is concave and $f(0)\ge0$,
\begin{equation*}
{\rm Tr\,} f(Z^*AZ)  \le {\rm Tr\,} Z^*f(A)Z.
\end{equation*} 

{\rm (2)} If $f$ is convex and $f(0)\le0$,
\begin{equation*}
{\rm Tr\,} f(Z^*AZ)  \ge {\rm Tr\,} Z^*f(A)Z.
\end{equation*} }

\vskip 5pt \noindent
{\bf Example 3.6.} Here, contrary to the Brown-Kosaki trace inequalities (Corollary 2.3(1)), the assumption $A\ge0$ is essential. For instance, in the convex case, consider $f(t)=t_+$,
\begin{equation*}
A=\begin{pmatrix} 1& 0 \\0&-1 \end{pmatrix}, \quad{\rm and}\quad Z=\begin{pmatrix} 2& 1 \\1& 2 \end{pmatrix}.
\end{equation*}
Then, we have ${\rm Tr}\,f(Z^*AZ)=3<5 ={\rm Tr}\,Z^*f(A)Z$. Of course, the assumption $A\ge0$ is also essential in Lemma 3.4.

\vskip 10pt
We turn to the proof of Theorem 3.5.

\vskip 10pt \noindent
{\bf Proof.} Of course, assertions (1) and (2) are equivalent. Let us prove (2). Since $Z$ is expansive we may assume that $f(0)=0$. By a limit argument we may then assume that
$$
f(t)=\lambda t + \sum_{i=1}^m \alpha_i(t-\beta_i)_+
$$
for a real $\lambda$ and some nonnegative reals $\{\alpha_i\}_{i=1}^m$ and $\{\beta_i\}_{i=1}^m$.
The result then follows from the linearity of the trace and Lemma 3.4. \qquad $\Box$  

\vskip 10pt
In order to extend Theorem 3.5(2) to all unitarily invariant norms, i.e.\ those norms $\Vert\cdot\Vert$
such that $\Vert UXV\Vert=\Vert X\Vert$ for all operators $X$ and all unitaries $U$ and $V$, we need a simple lemma. A family of positive operators $\{A_i\}_{i=1}^m$ is said to be {\it monotone} if there exists a positive operator $Z$ and a family  of nondecreasing nonnegative functions $\{f_i\}_{i=1}^m$ such that
$f_i(Z)=A_i$, $i=1,\dots m$.

\vskip 10pt \noindent
{\bf Lemma 3.7.} {\it Let $\{A_i\}_{i=1}^m$ be a monotone family of positive operators and  let $\{U_i\}_{i=1}^m$ be a family of unitary  operators. Then, for every unitarily invariant norm $\Vert\cdot\Vert$, we have
\begin{equation*}
\Vert \sum_i U_iA_iU_i^*\Vert \le \Vert \sum_i A_i\Vert.
\end{equation*}
}

\vskip 5pt \noindent
{\bf Proof.} By the Ky Fan dominance principle, it suffices to consider the Ky Fan $k$-norms $\Vert\cdot\Vert_{(k)}$ [1, pp.\ 92-3]. There exists a rank $k$ projection $E$ such that
$$
\Vert \sum_i U_iA_iU_i^*\Vert_{(k)} = \sum_i {\rm Tr\,} U_iA_iU_i^*E \le \sum_i \Vert A_i\Vert_{(k)}
=\Vert \sum_i A_i \Vert_{(k)}
$$
where the inequality comes from the maximal characterization of the Ky Fan norms and the last equality from the monotony of the family $\{A_i\}$. \qquad $\Box$

\vskip 10pt \noindent
{\bf Proposition 3.8.} {\it Let $A$ be a positive operator and $Z$ be an expansive operator. Assume that
$f$ is a nonnegative convex function defined on $[0,b]$, $b\ge \Vert Z^*AZ\Vert_{\infty}$. Assume also that $f(0)=0$. Then, for every unitarily invariant norm $\Vert\cdot\Vert$,
\begin{equation*}
\Vert f(Z^*AZ) \Vert \ge \Vert Z^*f(A)Z\Vert.
\end{equation*}
}

\vskip 5pt \noindent
{\bf Proof.} It suffices to consider the case when
$$
f(t)=\lambda t + \sum_{i=1}^m \alpha_i(t-\beta_i)_+
$$
for some nonnegative reals $\lambda$, $\{\alpha_i\}_{i=1}^m$ and $\{\beta_i\}_{i=1}^m$. By Lemma 3.4, we have
\begin{align*}
Z^*f(A)Z &= \lambda Z^*AZ + \sum_i Z^*\alpha_i(A-\beta_i I)_+Z \\
 &\le \lambda Z^*AZ + \sum_i U_i\alpha_i(Z^*AZ-\beta_i I)_+U_i^*
\end{align*} 
for some unitary operators $\{U_i\}_{i=1}^m$. Since $\lambda Z^*AZ$ and $\{\alpha_i(Z^*AZ-\beta_i I)_+\}_{i=1}^m$ form a monotone family, Lemma 3.7 completes the proof. \qquad $\Box$ 

\vskip 10pt \noindent
{\bf Theorem 3.9.} {\it Let $A$ be a positive operator, let $Z$ be an expansive operator and let
$f:[0,\infty)\longrightarrow[0,\infty)$ be a nondecreasing concave function. Then,
\begin{equation*}
\Vert f(Z^*AZ)\Vert_{\infty}  \le \Vert Z^*f(A)Z\Vert_{\infty}.
\end{equation*} 
}

\vskip 10pt \noindent
{\bf Proof.} Here, we assume that we are in the finite dimensional setting.

 Since $Z$ is expansive we may assume $f(0)=0$. By a continuity argument we may assume that $f$ is onto. Let $g$ be the reciprocal function. Note that $g$ is convex  and $g(0)=0$. By Proposition 3.8,
$$
\Vert g(Z^*AZ)\Vert_{\infty} \ge \Vert Z^*g(A)Z\Vert_{\infty}.
$$
Hence
$$
f(\Vert g(Z^*AZ)\Vert_{\infty}) \ge f(\Vert Z^*g(A)Z\Vert_{\infty}).
$$
Equivalently,
$$
\Vert Z^*AZ\Vert_{\infty} \ge \Vert f(Z^*g(A)Z)\Vert_{\infty},
$$
so, letting $B=g(A)$,
$$
\Vert Z^*f(B)Z\Vert_{\infty} \ge \Vert f(Z^*BZ)\Vert_{\infty},
$$
proving the result because $A\longrightarrow g(A)$ is onto. \qquad $\Box$

\vskip 10pt 
Our next result is a straightforward application of Theorem 2.2.

\vskip 10pt \noindent
{\bf Proposition 3.10.} {\it Let $A$ be a positive operator and $Z$ be an expansive operator. Assume that
$f$ is a nonnegative function defined on $[0,b]$, $b\ge \Vert Z^*AZ\Vert_{\infty}$. Then:
\vskip 5pt
 $(1)$ If $f$ is concave nondecreasing,
\begin{equation*}
\det f(Z^*AZ)  \le \det Z^*f(A)Z.
\end{equation*} 
\vskip 5pt
 $(2)$ If $f$ is convex increasing and $f(0)=0$,
\begin{equation*}
\det f(Z^*AZ)  \ge \det Z^*f(A)Z.
\end{equation*}
}

\vskip 5pt \noindent
{\bf Proof.} For instance, consider the concave case. By Theorem 2.2, there exists a unitary operator $U$ such that $Z^{*-1}f(Z^*AZ)Z^{-1}\le Uf(A)U^*$; hence the result follows. \qquad $\Box$

\vskip 10pt 
We note the following fact about operator convex functions:

\vskip 10pt \noindent
{\bf Proposition 3.11.} {\it Let $f:[0,\infty)\longrightarrow[0,\infty)$ be a one to one continuous function with $f(0)=0$ and $f(\infty)=\infty$. The following statements are equivalent:
\vskip 5pt
{\rm \ (i)} f(t) is operator convex.
\vskip 5pt
{\rm  (ii)} 1/f(1/t) is operator convex.
}

\vskip 10pt \noindent
{\bf Proof.} Since the map $f(t)\longrightarrow\Psi(f)(t)=1/f(1/t)$ is an involution on the set of all one to one continuous functions $f$ on $[0,\infty)$ with $f(0)=0$ and $f(\infty)=\infty$, it suffices to check that ${\rm (i)}\Rightarrow{\rm (ii)}$. But, by the Hansen-Pedersen inequality [6], (i) is equivalent to
\begin{equation*}
f(Z^*AZ)\le Z^*f(A)Z \tag{9}
\end{equation*}
for all $A\ge0$ and all contractions $Z$. By a limit argument, it suffices to require (9) when both $A$
and $Z$ are invertible. Then, as $t\longrightarrow t^{-1}$ is operator decreasing, (9) can be written
\begin{equation*}
f^{-1}(Z^*AZ)\ge Z^{-1}f^{-1}(A)Z^{*-1}, 
\end{equation*}
or
\begin{equation*}
f^{-1}(A)\le Zf^{-1}(Z^*AZ)Z^*, 
\end{equation*}
but  this is the same as saying that (9) holds for $\Psi(f)$, therefore  $\Psi(f)$ is operator convex. \qquad $\Box$

\vskip 10pt
We wish to sketch another proof of Proposition 3.10.  By a result of Hansen and Pedersen [6], for a continuous function $f$ on $[0,\infty)$, the following conditions are equivalent:
\vskip 5pt
\ (i) $f(0)\le 0$ and $f$ is operator convex.
\vskip 5pt
(ii) $t\longrightarrow f(t)/t$ is operator monotone on $(0,\infty)$.
\vskip 5pt\noindent
Using the operator monotony of $t\longrightarrow 1/t$ on $(0,\infty)$, we note that if $f(t)$ satisfies to (ii), then so does $1/f(1/t)$. This proves Proposition 3.10.

\vskip 10pt \noindent
{\bf Question 3.12.} Does Theorem 3.9 extend to all nonnegative concave functions on $[0,b]$ and/or to all unitarily invariant norms ?

\vskip 20pt \noindent
{\large\bf 4. Addendum}

\vskip 10pt 
There exist several inequalities involving $f(A+B)$ and $f(A)+f(B)$ where $A$, $B$ are Hermitians and $f$ is a function with special properties. We wish to state and prove one of the most basic results in this direction which can be derived from a more general result due to Rotfel'd (see [1, p.\ 97]). The simple proof given here is inspired by that of Theorem 3.5.

\vskip 10pt \noindent
{\bf Proposition 4.1.} {\rm (Rotfel'd)} {\it Let $A$, $B$ be positive operators.
 \vskip 5pt
 {\rm(1)} If $f$ is a convex  nonnegative function on $[0,\infty)$ with $f(0)\le 0$, then 
\begin{equation*}
{\rm Tr}\, f(A+B) \ge {\rm Tr}\,f(A)+ {\rm Tr}\,f(B).
\end{equation*}
\vskip 5pt
 {\rm(2)} If $f$ is a concave nonnegative function on $[0,\infty)$, then 
\begin{equation*}
{\rm Tr}\, g(A+B) \le {\rm Tr}\, g(A)+ {\rm Tr}\,g(B).
\end{equation*}
}

\vskip 5pt \noindent
{\bf Proof.} By limit arguments, we may assume that we are in the finite dimensional setting. Since, on any compact interval $[a,b]$, $a>0$, we may write $g(x)=\lambda x-f(x)+\mu$ for some scalar $\lambda,\mu\ge 0$ and some convex function $f$ with $f(0)=0$, it suffices to consider the convex case. Clearly we may assume $f(0)=0$. Then, $f$ can be uniformly approximated, on any compact interval, by
a positive combination of functions $f_{\alpha}(x)=\max\{0,x-\alpha\}$, $\alpha>0$.

Therefore, still using the notation $S_+$ for the positive part of the Hermitian operator $S$, we need only to show that
$$
{\rm Tr}\,(A+B-\alpha)_+ \ge {\rm Tr}\,(A-\alpha)_+ + {\rm Tr}\,(B-\alpha)_+ .
$$
To this end, consider an orthonormal basis $\{e_i\}_{i=1}^n$ of eigenvectors for $A+B$. We note that:

(a) If $\langle e_i,(A+B-\alpha)_+e_i\rangle=0$, then $\langle e_i,(A+B-\alpha)_+e_i\rangle\le\alpha$ so that we also have $\langle e_i,(A-\alpha)_+e_i\rangle=\langle e_i,(B-\alpha)_+e_i\rangle=0$.

(b) If $\langle e_i,(A+B-\alpha)_+e_i\rangle>0$, then we may write
$$
\langle e_i,(A+B-\alpha)_+e_i\rangle=\langle e_i,Ae_i\rangle-\theta\alpha +\langle e_i,Be_i\rangle -(1-\theta)\alpha
$$
for some $0\le \theta\le 1$ chosen in such a way that $\langle e_i,Ae_i\rangle-\theta\alpha \ge0$ and
$\langle e_i,Be_i\rangle -(1-\theta)\alpha\ge0$. Hence, we have
$$
\langle e_i,Ae_i\rangle-\theta\alpha = \langle e_i,(A-\theta\alpha)_+e_i\rangle \ge
\langle e_i,(A-\alpha)_+e_i\rangle 
$$
and
$$
\langle e_i,Be_i\rangle-(1-\theta)\alpha = \langle e_i,(B-(1-\theta)\alpha)_+e_i\rangle \ge
\langle e_i,(B-\alpha)_+e_i\rangle 
$$
by using the simple fact that for two commuting Hermitian operators $S,T$, $S\le T\Rightarrow S_+\le T_+$.

From (a) and (b) we derive the desired trace inequality by summing over $i=1,\dots n$. \qquad $\Box$

\vskip 5pt
{\bf References}
\vskip 5pt
\noindent
{\small 
\noindent
[1] R. Bhatia, Matrix Analysis, Springer, Germany, 1996

\noindent
[2] J.-C. Bourin, Total dilation II, Linear Algebra Appl. 374/C  (2003) 
19-29.

\noindent
[3] L.\ G.\ Brown and  H.\ Kosaki, Jensen's inequality is semi-finite von Neumann algebras, J. Operator Theory 23 (1990) 3-19.

\noindent
[4] C.\ Davis, A Schwarz inequality for convex operator functions, Proc.\ Amer.\ Math.\ Soc.\  8 (1957) 42-44.

\noindent
[5] F.\ Hansen, An operator inequality, Math. Ann. 258 (1980) 249-250.

\noindent
[6] F.\ Hansen and G.\ K.\ Pedersen, Jensen's inequality for operator sand Lowner's Theorem, Math. Ann. 258 (1982) 229-241.

\noindent
[7] F.\ Hansen and G.\ K.\ Pedersen, Jensen's operator inequality, Bull.\ London Math.\ Soc.\ 35 (2003) 553-564.

\noindent
[8] B.\ Simon, Trace Ideals and Their Applications LMS lecture note, 35 Cambridge Univ.\ Press, Cambridge,
1979

\noindent
[9] X.\ Zhan, Matrix Inequalities, LNM 1790, Springer, Berlin, 2002.
}

\end{document}